\documentstyle[12pt]{article}
\begin{document}
\newpage
\pagestyle{empty}
\setcounter{page}{0}
\renewcommand{\thesection}{\Roman{section}}
\renewcommand{\theequation}{\thesection.\arabic{equation}}
\newcommand{\sect}[1]{\setcounter{equation}{0}\section{#1}}
\newfont{\twelvemsb}{msbm10 scaled\magstep1}
\newfont{\eightmsb}{msbm8}
\newfont{\sixmsb}{msbm6}
\newfam\msbfam
\textfont\msbfam=\twelvemsb
\scriptfont\msbfam=\eightmsb
\scriptscriptfont\msbfam=\sixmsb
\catcode`\@=11
\def\Bbb{\ifmmode\let\next\Bbb@\else
  \def\next{\errmessage{Use \string\Bbb\space only in math mode}}\fi\next}
\def\Bbb@#1{{\Bbb@@{#1}}}
\def\Bbb@@#1{\fam\msbfam#1}
\newfont{\twelvegoth}{eufm10 scaled\magstep1}
\newfont{\tengoth}{eufm10}
\newfont{\eightgoth}{eufm8}
\newfont{\sixgoth}{eufm6}
\newfam\gothfam
\textfont\gothfam=\twelvegoth
\scriptfont\gothfam=\eightgoth
\scriptscriptfont\gothfam=\sixgoth
\def\frak{\frak@}
\def\frak@#1{{\fam\gothfam{{#1}}}}
\def\frak@@#1{\fam\gothfam#1}
\catcode`@=12
%
%
%
\def\CC{{\Bbb C}}
\def\NN{{\Bbb N}}
\def\QQ{{\Bbb Q}}
\def\RR{{\Bbb R}}
\def\ZZ{{\Bbb Z}}
\def\cA{{\cal A}}          \def\cB{{\cal B}}          \def\cC{{\cal C}}
\def\cD{{\cal D}}          \def\cE{{\cal E}}          \def\cF{{\cal F}}
\def\cG{{\cal G}}          \def\cH{{\cal H}}          \def\cI{{\cal I}}
\def\cJ{{\cal J}}          \def\cK{{\cal K}}          \def\cL{{\cal L}} 
\def\cM{{\cal M}}          \def\cN{{\cal N}}          \def\cO{{\cal O}}
\def\cP{{\cal P}}          \def\cQ{{\cal Q}}          \def\cR{{\cal R}} 
\def\cS{{\cal S}}          \def\cT{{\cal T}}          \def\cU{{\cal U}}
\def\cV{{\cal V}}          \def\cW{{\cal W}}          \def\cX{{\cal X}}
\def\cY{{\cal Y}}          \def\cZ{{\cal Z}}
\def\qed{\hfill \rule{5pt}{5pt}}
\def\id{\mbox{id}}
\def\ggo{{\frak g}_{\bar 0}}
\def\uqggo{\cU_q({\frak g}_{\bar 0})}
\def\uqggp{\cU_q({\frak g}_+)}
\def\typeA{{\em type $\cA$}}
\def\typeB{{\em type $\cB$}}
\newtheorem{lemma}{Lemma}
\newtheorem{prop}{Proposition}
\newtheorem{theo}{Theorem}
\newtheorem{Defi}{Definition}
$$
\;
$$
$$
\;
$$
$$
\;
$$
$$
\;
$$
\vfill
\vfill
\begin{center}

{\LARGE {\bf {\sf 
The Twisted Heisenberg Algebra ${\cal U}_{h,w}({\cal H}(4))$
}}} \\[2cm]

{\large Boucif Abdesselam\footnote{boucif@orphee.polytechnique.fr}$^{}$}

\smallskip 
\smallskip 
\smallskip 
\smallskip 

{\em  \footnote{Laboratoire Propre du CNRS UPR A.0014}$^{}$Centre de 
Physique Th\'eorique, Ecole Polytechnique,\\ 91128 Palaiseau Cedex, France.}

\smallskip 
\smallskip 
\smallskip 
\smallskip 
\smallskip 
\smallskip 
\end{center}

\vfill

\begin{abstract}
\noindent A two parametric deformation of the enveloping Heisenberg algebra 
${\cal H}(4)$
which appear as a combination of the standard and a nonstandard quantization 
given by Ballesteros and Herranz is defined and proved to be Ribbon Hopf 
algebra. The universal ${\cal R}$-matrix and its associated quantum group are 
constructed. New solution of Braid group are obtained. The contribution of 
these parameters in invariants of links and WZW model are analyzed. General 
results for twisted Ribbon Hopf algebra are derived. 
\end{abstract}

\vfill
\vfill

\rightline{CPTH-S464.1196}

$$
\;
$$
$$
\;
$$

\newpage 

\pagestyle{plain}

\sect{Introduction}

An enveloping Lie algebra ${\cal U}({\cal G})$ has many quantizations ($\geq 
2$): The first one is called the Drinfeld-Jimbo quantization \cite{Drin,Jim}, 
whereas the other ones are called the nonstandard quantizations 
\cite{Manin}\footnote{In the case of ${\cal U}(sl(2))$, the nonstandard 
quantization \cite{Ohn} was obtained as contraction of the Drinfeld-Jimbo 
one (see \cite{Iran}). The universal ${\cal R}$-matrix of the nonstandard 
algebra ${\cal U}_{h}(sl(2))$ was obtained \cite{Vlad,Shar}.}. 
Recently there is much interest in studies relating to various aspects of 
the nonstandard quantizations. 

\smallskip \smallskip \smallskip

The standard $q$-Heisenberg algebras and their universal 
${\cal R}$-matrices \cite{Cel} has recently attracted wide attention. The 
use of $q$-Heisenberg algebras to describe composite particles \cite{Avan}, 
the description of certain classes of exactly solvable potentials in terms of 
a $q$-Heisenberg dynamical symmetry \cite{Spi}, the link between deformed 
oscillator algebras and superintegrable systems \cite{Bon2,Bon1} and the 
relations
between these deformed algebras and $q$-orthogonal polynomials \cite{Vinet} 
are the most attractive examples. In reference \cite{Gom}, the quasitriangular
$q$-oscillator algebra has been found to be related to Yang-Baxter systems 
and link invariants. 
         
The Heisenberg algebras has also many nonstandard quantizations. Recently, the 
coboundary Lie bialgebras and their corresponding Poisson-Lie structures has 
been constructed for the Heisenberg algebra. The quantum nonstandard 
Heisenberg algebras are derived from these bialgebras by using the Lyakhovsky 
and Mudrov formalism and for some cases, quantizations at both algebra an 
group levels have been obtained, including their universal ${\cal R}$-matrices
\cite{Her}. 
     
\smallskip \smallskip \smallskip

Let us mention that the combination of some nonstandard quantizations to the 
standard ones permits to equip the enveloping algebras with multiparameter 
quantizations. Following this way, we define a new two parametric quantization 
of the Heisenberg algebras ${\cal U}_{h,w}({\cal H}(4))$ by combination of the 
standard and a nonstandard quantization given in \cite{Her}. This type of 
deformed bosons can be expected to build up $[h,w]$-boson realizations of a 
possible two parametric ${\cal U}_{h,w}(sl(2))$ algebra which not exit in
the literature at our knowlodge. 

\smallskip \smallskip \smallskip

 One of the instances in which both the parameters of the nonstandard 
Heisenberg algebra become relevant is in the context of reflection 
equation 
\begin{equation}
R_{12}K_1R_{21}K_2 = K_2R_{12}K_1R_{21}\,.
\label{re}
\end{equation}
It was introduced in~\cite{C} in the study of two particle 
scattering on a half-line, where the matrix $K$ described the 
reflection of a particle at the end point. It is also known~\cite{AFS}
to have other applications, such as, a generalization of the inverse 
scattering method to the case of non-ultralocal commutation 
relations and a lattice regulatized version of Kac-Moody algebras. 
The reflection algebra is closely connected to the quantum group and 
the matrix $K$ may be realized~\cite{AF} from the knowledge of the 
Lax operators $L^{\pm}$~:
\begin{equation}
K = S\left( L^- \right) L^+\,.
\label{k}
\end{equation} 
One viewpoint advocated by Majid~\cite{M} is the `transmutation 
procedure' converting the quantum groups to the braided groups, 
which may be looked simply as a generalization of the supergroups 
with $\pm$ Bose-Fermi statistics being replaced by braid 
statistics.  There are indications~\cite{FRS} that particleds of braid 
statistics arise in low dimensional field theory.  As the 
construction of the above braided algebra can be done for any 
regular invertible $R$-matrix, it is particularly relevant to study 
the deformed Heisenberg algebras, which are closely connected to 
the statistics problem.  

From the knowledge of the universal $R$-matrix of 
$U_{h,\omega}(H(4))$, we can construct the Lax operators $L^{\pm}$, 
which, as a consequence of~(\ref{k}), generates the reflection 
operator $K$.  This necessarily depends on both deformation 
parameters $(h, \omega )$ and, in the limit $\omega \longrightarrow 
0$ limit, yields the $K$ operator of the $U_h(H(4))$ algebra.  In 
the present context we would just like to point out that in the 
braid statistics problem, which appears for a chain of $(h, \omega 
)$-deformed oscillators, both deformation parameters are physically 
important.  

\smallskip \smallskip \smallskip

The algebra ${\cal U}_{h,w}({\cal H}(4))$ can also be interpreted 
as two parametric deformation of an extended $(1+1)$-poincar\'e algebra. 
Others interesting applications related to braid groups and special 
functions theory of our two parametric quantization can be found.   
 
\smallskip \smallskip \smallskip

The main purpose of this work is fourfold: 

\smallskip

{\bf (i)} To provide the enveloping 
Heisenberg algebras with a two parametric deformation which appear as a 
combination of the standard quantization and a nonstandard one \cite{Her}. 

\smallskip

{\bf (ii)} We derive its universal ${\cal R}^{h,w}$-matrix, we use it to make 
explicit connection with the formalism of matrix quantum pseudogroups due to 
Woronowicz \cite{Wor1,Wor2} and to define the Hopf algebra of the 
representative elements. 

\smallskip

{\bf (iii)} Using the results obtained in \cite{Gom},
this two parametric deformation is proved to be a ribbon algebra. General 
results for a twisted ribbon Hopf algebra are also derived. 

\smallskip

{\bf (iv)} We give
a new solution of the braid group ${\cal B}_{m}$ and we analyze briefly the 
contribution of these two parameters in WZW model and in invariants of links.
         
\smallskip \smallskip \smallskip

This paper is organized as follows: In section II, we give the definitions of 
the standard and nonstandard quantizations of the Heisenberg algebra, their 
central elements and their universal ${\cal R}$-matrices. In section III, 
we introduce the two parametric deformation and its quasitriangular 
Hopf structure. New solution of the braid group ${\cal B}_{m}$ and the 
contribution of these parameters in the theory of invariants of links is also 
discussed. In the section IV, we determine the infinitesimal generators of the 
matrix pseudogroup and we calculate both commutations relations, coproducts, 
counits and antipodes. The two parametric deformation is proved to be a ribbon 
Hopf algebra in section V. General results are also established. Some comments 
concerning the WZW model are also presented. Finally, we conclude with some 
remarks and perspectives in section VI.    
         
\sect{Standard and Nonstandard Heisenberg Algebras and their 
Universal ${\cal R}$-matrices} 

In this paper, $h$ and $w$ are arbitrary complex numbers. We denote by $n$, 
$e$, $a^{+}$ and $a$ the generators of the Heisenberg algebra ${\cal H}(4)$
that satisfy the commutation relations  
\begin{eqnarray}   
&& [a,a^{+}]=e,\qquad\qquad [n,a]=-a,
\qquad\qquad [n,a^{+}]=a^{+}, \qquad\qquad [e,\;.\;]=0 \qquad\qquad 
\end{eqnarray}
and denote by ${\cal U}({\cal H}(4))$ the enveloping algebra of ${\cal H}(4)$. 
The Heisenberg group is denoted by $H(4)$. Let 
\begin{eqnarray}   
r=a\otimes a^{+}-e\otimes n
\end{eqnarray}
and introduce the notation $r_{ij}$, $1\leq i < j\leq  3$, where 
$r_{12}=\sum_{i} a_{i}\otimes b_{i}\otimes 1$, $r_{23}=\sum_{i} 1
\otimes a_{i}\otimes b_{i}$ and $r_{13}=\sum_{i} a_{i}\otimes 1\otimes b_{i}$, 
if $r=\sum_{i} a_{i}\otimes b_{i}\in{\cal U}({\cal H}(4))\otimes 
{\cal U}({\cal H}(4))$. The element (II.2) solves the classical 
Yang-Baxter equation (CYBE)
\begin{eqnarray}  
[r_{12}, r_{13}]+[r_{12}, r_{23}]+ [r_{13}, r_{23}]=0  
\end{eqnarray}
and its called a classical $r$-matrix. The quantum Hopf algebra 
${\cal U}_{h}({\cal H}(4))$ which quantizes the standard bialgebra generated 
by the classical $r$-matrix (II.2) is defined as: 

\begin{Defi} 
The quantum standard algebra ${\cal U}_{h}({\cal H}(4))$ is the 
unital associative algebra with the generators $N$, $E$, $A^{+}$, $A$ and 
the relations
\begin{eqnarray}
&& [A, A^{+}]=\displaystyle{\sinh(hE)\over h},
\qquad [N, A^{+}]= A^{+},\qquad [N,A]= -A,\qquad [E,\;.\;]= 0. \qquad 
\end{eqnarray}
The algebra ${\cal U}_{h}({\cal H}(4))$ admits a Hopf structure with 
coproducts, counits and antipodes determined by
\begin{equation}
\begin{array}{lll}
\bigtriangleup_{h}(N)=N\otimes 1 + 1\otimes N,&  \qquad S_{h}(N)=-N,& \qquad 
\varepsilon_{h}(N)=0,\nonumber  \\
\bigtriangleup_{h}(E)=E\otimes 1 + 1 \otimes E,& \qquad S_{h}(E)=-E,& \qquad 
\varepsilon_{h}(E)=0,\nonumber \\
\bigtriangleup_{h}(A^{+})=A^{+} \otimes 1 + 1 \otimes A^{+},
& \qquad S_{h}(A^{+})=-A^{+}, &\qquad \varepsilon_{h}(A^{+})=0,
\nonumber  \\
\bigtriangleup_{h}(A)=A \otimes e^{hE}+e^{-hE}\otimes A,&
\qquad S_{h}(A)=-A,& \qquad 
\varepsilon_{h}(A)=0.
\end{array}
\end{equation}
In this point of view, the generators $N$, $E$ and $A^{+}$ are primitive 
elements.
\end{Defi}

There exist another element belonging to the center of the algebra 
${\cal U}_{h}({\cal H}(4))$ given by     
\begin{eqnarray}
C_{h}= N\;{\sinh{ h E} \over h} - {1\over 2}(A^{+}A + 
A\;A^{+}).  
\end{eqnarray}
From the relations (II.5), we have $S_{h}^{2}=\id$. The quantum Hopf algebra 
(II.4) and (II.5) is equivalent to the structure defined by Celeghini et all 
\cite{Cel}. The coproduct $\bigtriangleup_{h}$, counit $\varepsilon_{h}$ and 
antipode $S_{h}$ (II.5) are related to Celeghini et all coproduct 
$\bigtriangleup_{c}$, counit $\varepsilon_{c}$ and antipode $S_{c}$ by  
\begin{eqnarray}  
&& \bigtriangleup_{h} =
{\cal A}^{-1}\bigtriangleup_{c}{\cal A},\qquad\qquad\varepsilon_{h}
=\varepsilon_{c},\qquad\qquad S_{h}=S_{c},  
\end{eqnarray}
where, 
${\cal A}=e^{-{h\over 2}(E\otimes N-N\otimes E)}$. The classical $r$-matrix 
associated to the quantum structure defined in \cite{Cel} is given by  
\begin{eqnarray}   
r= a\otimes a^{+} -{1 \over 2}(e\otimes n + n\otimes e).   
\end{eqnarray}     

Recall that a quasitriangular Hopf algebra is a pair $({\cal U}, {\cal R})$ 
where ${\cal U}$ is a Hopf algebra and ${\cal R}$ is invertible element 
obeying to the following axioms 
\begin{eqnarray}     
&& {\cal R}\bigtriangleup(u) {\cal R}^{-1}= \sigma \circ \bigtriangleup (u),
\qquad\qquad u\in {\cal U}, \\
&& (\bigtriangleup \otimes \id )({\cal R})={\cal R}_{13}{\cal R}_{23},
\qquad\qquad(\id\otimes\bigtriangleup)({\cal R})={\cal R}_{13}{\cal R}_{12}, 
\end{eqnarray}
where, if ${\cal R}=\sum_{i}a_i\otimes b_i\in{\cal U}\otimes{\cal U}$, 
we denote ${\cal R}_{12}=\sum_{i}a_i\otimes b_i\otimes 1\in {\cal U}
\otimes{\cal U}\otimes{\cal U}$, ${\cal R}_{13}= \sum_{i}a_i\otimes 1 
\otimes b_i$, ${\cal R}_{23}= \sum_{i}1\otimes a_i \otimes b_i$ and 
$\sigma$ is the flip operator $\sigma (a\otimes b)= b\otimes a$. The relation 
(II.9) indicates that ${\cal R}$ being a intertwining operator on the 
coproduct $\bigtriangleup$. ${\cal R}$ is called an universal 
${\cal R}$-matrix and satisfies the quantum Yang-Baxter equation (QYBE)
\begin{eqnarray}   
{\cal R}_{12}{\cal R}_{13}{\cal R}_{23}=
{\cal R}_{23}{\cal R}_{13}{\cal R}_{12}.
\end{eqnarray}

Now, given a quasitriangular Hopf algebra $({\cal U},{\cal R})$ and an 
invertible element ${\cal F}\in {\cal U}\otimes{\cal U}$ satisfying 
the following conditions:  
\begin{eqnarray}
&&(\bigtriangleup_{h}\otimes\id)({\cal F}){\cal F}_{12} =
(\id \otimes \bigtriangleup_{h})({\cal F}){\cal F}_{23} ,  \\
&& {\cal F}_{12}^{-1}(\bigtriangleup_{h}\otimes \id )({\cal F}^{-1})=
{\cal F}_{23}^{-1}(\id \otimes \bigtriangleup_{h})({\cal F}^{-1}), \\
&&(\varepsilon_{h}\otimes\id)({\cal F})=(\id\otimes\varepsilon_{h})
({\cal F}),\\
&&(\varepsilon_{h}\otimes\id)({\cal F}^{-1})=(\id\otimes\varepsilon_{h})
({\cal F}^{-1})
\end{eqnarray}
one can form a new quasitriangular Hopf algebra ${\cal U}_{\cal F}$ by 
twisting ${\cal U}$: \cite{Drinfeld} ${\cal U}_{\cal F}$ retains the vector 
space and multiplication of ${\cal U}$ while its quasitriangular Hopf 
structure is given by
\begin{eqnarray} 
&& \bigtriangleup_{\cal F}\equiv {\cal F}^{-1}\bigtriangleup {\cal F}, 
\nonumber \\
&& \varepsilon_{\cal F}= \varepsilon,\nonumber \\
&& S_{\cal F}= v^{-1}Sv, \nonumber \\
&& {\cal R}^{\cal F}= {\cal F}_{21}^{-1}{\cal R}{\cal F}  \qquad\qquad 
\end{eqnarray}
where, the element $v$ and its inverse are obtained from the invertible 
element ${\cal F}\in {\cal U}\otimes {\cal U}$ as follows:   
\begin{eqnarray}
&& v=m(S\otimes \id)({\cal F}) \qquad\qquad 
v^{-1}=m(\id \otimes S)({\cal F}^{-1}). 
\end{eqnarray}
If $S^{2}=\id$, the twisted quasitriangular Hopf algebra ${\cal U}_{\cal F}$ 
thus not conserve the same property, namely, 
\begin{eqnarray}
S_{\cal F}^{2}(a)=v^{-1}S(v)a S(v^{-1})v.
\end{eqnarray} 

\begin{prop} Let $({\cal U},{\cal R})$ a quasitriangular Hopf algebra and 
an element $ {\cal F}\in {\cal U}\otimes {\cal U}$ satisfying 
the conditions (II.12)-(II.15). The element $v$ defined in (II.17) has the 
following properties: 
\begin{eqnarray}
&&\varepsilon(v)=1, \qquad\qquad
\bigtriangleup(v)=(S\otimes S)({\cal F}_{21}^{-1})\;(v\otimes v)\; 
{\cal F}^{-1},   
\end{eqnarray}
and, for the element $v^{-1}S(v)$, we have
\begin{eqnarray}
\bigtriangleup(v^{-1}S(v)) = {\cal F} (v^{-1}S(v)\otimes v^{-1}S(v))
(S^2 \otimes S^2)({\cal F}^{-1}).
\end{eqnarray}
If $S^2=\id$, the relation (II.20) reads
\begin{eqnarray}
\bigtriangleup_{\cal F}(v^{-1}S(v)) = v^{-1}S(v)\otimes v^{-1}S(v).
\end{eqnarray}
\end{prop}

\noindent {\bf Proof:} The equations (II.19) arise by direct calculations.
The relations (II.20) and (II.21) are derived using (II.19). \qed

\begin{prop} The Hopf algebra ${\cal U}_{h}({\cal H}(4))$ is quasitriangular. 
The universal ${\cal R}$-matrix has the following form 
\begin{eqnarray}   
&&\displaystyle{\cal R}^{h}=\exp(-2h\;E\otimes N)\;\exp\biggl(2h\;e^{hE}A
\otimes A^{+}\biggr)
\end{eqnarray}
and satisfies the Quantum 
Yang-Baxter equation (II.11).    
\end{prop}
\noindent {\bf Proof:} Similar to the proof given in \cite{Cel} 
(see proposition II.9). The ${\cal R}$-matrix (II.22) is related to 
${\cal R}_{c}$-matrix given in \cite{Cel} by the following twist 
\begin{eqnarray}  
{\cal R}^{h}=A_{21}^{-1}\;{\cal R}^{c}\;{\cal A}. 
\end{eqnarray}
\qed

Now, let us consider the nonstandard classical $r$-matrix 
\begin{eqnarray}   
r= n\otimes a^{+} - a^{+}\otimes n,  
\end{eqnarray}
which solves the classical Yang-Baxter equation (II.3). The quantum Hopf 
algebra which quantizes the nonstandard bialgebra generated by the 
classical $r$-matrix (II.24) was given by Ballesteros and Herranz \cite{Her}: 
 
\begin{Defi} 
The quantum nonstandard algebra ${\cal U}_{w}({\cal H}(4))$ is the unital 
associative algebra generated by $N$, $E$, $A^{+}$ and $A$, satisfying the 
commutations relations
\begin{eqnarray}
  && [A,A^{+}]=E\;e^{wA^{+}},\qquad 
[N,A^{+}]= {e^{wA^{+}}-1\over w},\qquad
 [N,A]= - A, \qquad [E,\;.\;]=0.\qquad 
\end{eqnarray}
The algebra ${\cal U}_{w}({\cal H}(4))$ has the following Hopf structure
\begin{equation}
\begin{array}{lll}
\bigtriangleup_{w}(N)=N\otimes e^{wA^{+}}  + 1\otimes N, & 
S_{w}(N)=- N e^{-w A^{+}}, &  \varepsilon_{w} (N)= 0, \nonumber  \\
\bigtriangleup_{w} (E)=E\otimes 1 + 1 \otimes E, & 
S_{w}(E)=- E, &  \varepsilon_{w} (E)= 0,   \nonumber  \\
\bigtriangleup_{w} (A^{+})=A^{+} \otimes 1 + 1 \otimes A^{+},
& S_{w}(A^{+})=- A^{+}, & \varepsilon_{w} (A^{+})= 0
\end{array}
\end{equation}
and
\begin{eqnarray}
&& \bigtriangleup_{w}(A)=A \otimes e^{w A^{+}}+1\otimes A+w N \otimes 
E e^{w A^{+}},\qquad S_{w}(A)=-Ae^{-w A^{+}}+wNEe^{-w A^{+}}, 
\nonumber \\
&& \qquad\qquad\qquad\qquad\qquad\qquad\qquad\qquad \varepsilon_{w}(A)=0.
\phantom{a \over b}
\end{eqnarray} 
The generators $E$ and $A_{+}$ are primitives.  
\end{Defi}

We will not consider here the others nonstandard deformations. The Casimir 
element of ${\cal U}_{w}({\cal H}(4))$ is given by  \cite{Her}
\begin{eqnarray}
C_{w}= N\;E + {e^{-w A_{+}} -1 \over 2\;w} A+ 
A {e^{-w A_{+}} -1 \over 2\;w}.  
\end{eqnarray} 

\begin{prop}
The Hopf algebra ${\cal U}_{w}({\cal H}(4))$ is quasitriangular. The universal
${\cal R}$-matrix has the following form \cite{Her}: 
\begin{eqnarray}   
\displaystyle {\cal R}^{w}= \exp(-w\;A^{+}\otimes N)\;\exp (w\;N\otimes 
A^{+})
\end{eqnarray}
and satisfies the quantum Yang-Baxter equation (II.11).
\end{prop}

In passing, let us mention that ${\cal R}^{w}= {\cal F}_{21}^{-1}
{\cal F}$, with ${\cal F}=\exp(wN\otimes A^{+})$ and $\bigtriangleup_{w}
={\cal F}^{-1}\bigtriangleup_{0}{\cal F}$, where $\bigtriangleup_{0}$ 
is the coproduct of the enveloping algebra ${\cal U}({\cal H}(4))$, 
extended to ${\cal U}({\cal H}(4))[[w]]$, namely 
\begin{eqnarray}
&& \bigtriangleup_{0}(N)=N\otimes 1 + 1\otimes N,\nonumber  \\
&& \bigtriangleup_{0}(E)=E\otimes 1 + 1\otimes E,\nonumber \\
&& \displaystyle \bigtriangleup_{0}\biggl({1-e^{-wA^{+}}\over w}\biggr)=
\biggl({1-e^{-wA^{+}}\over w}\biggr)\otimes 1 + 1\otimes 
\biggl({1-e^{-wA^{+}}\over w}\biggr),\nonumber \\
&& \bigtriangleup_{0}(A)=A \otimes 1+1\otimes A.
\end{eqnarray}
${\cal F}$ appear as an element which deforms the coproduct 
$\bigtriangleup_{0}$ of ${\cal U}({\cal H}(4))$ to the coproduct 
$\bigtriangleup_{w}$ of ${\cal U}_{w}({\cal H}(4))$. Furthermore, the 
antipodes and counits are
\begin{eqnarray}
&&S_{0}(N)=-N,\qquad\qquad S_{0}(A)=-A,\qquad\qquad S_{0}\biggl({1-e^{-wA^{+}}
\over w}\biggr)=-\biggl({1-e^{-wA^{+}}\over w}\biggr), \nonumber \\ 
&& \qquad\qquad\qquad\qquad  \varepsilon_{0}(N)= \varepsilon_{0}(A)= 
\varepsilon_{0}\biggl({1-e^{-wA^{+}}\over w}\biggr)=0. 
\end{eqnarray}

Let us just mention, that there exist others solutions of the classical 
Yang-Baxter equation given by 
\begin{eqnarray}
r(\mu,\nu)=\mu(a\otimes e-e\otimes a)+\nu(a^{+}\otimes e-e\otimes a^{+}),
\qquad\qquad \mu,\nu \in \CC.       
\end{eqnarray} 
The quantum Hopf algebra which quantizes the nonstandard bialgebra generated 
by the classical $r$-matrix $r(\mu,0)$ can be obtained as contraction 
limit from ${\cal U}_{h}(sl(2))\otimes u(1)$. Another 
interesting solution will be the subject of the following section.

\sect{A Two Parametric Deformation of ${\cal U}({\cal H}(4))$ and Links 
Invariants}

In this section, we analyse the algebra which corresponds to a combination
of the standard deformation (II.4) and the nonstandard one (II.25). Let us 
consider the element
\begin{eqnarray} 
r(h,w)=2h(a\otimes a^{+}-e\otimes n)+w(n\otimes a^{+}-a^{+}\otimes n), 
\qquad h,w \in \CC
\end{eqnarray}
which satisfies the classical Yang-Baxter equation (II.3). In this case, 
the quantum Hopf algebra ${\cal U}_{h,w}({\cal H}(4))$ which quantizes the 
bialgebra generated by the $r(h,w)$-matrix (III.1) will be  characterized by 
two parameters $h$ and $w$ associated respectively to the standard $r$-matrix 
$r(h,0)=2h (a\otimes a^{+}-e\otimes n)$ and to the nonstandard $r$-matrix 
$r(0,w)=w(n\otimes a^{+}-a^{+}\otimes n)$:      

\begin{prop}
The two parametric deformed algebra ${\cal U}_{h,w}({\cal H}(4))$ is an 
associative algebra over $\CC$ generated by $N$, $E$, $A^{+}$ and $A$, 
satisfying the commutations relations
\begin{eqnarray}
  && [A,A^{+}]={\sinh (hE) \;e^{w A^{+}}\over h}, 
\qquad [N,A^{+}]= {e^{wA^{+}}-1\over w},\qquad  [N,A]= - A 
\end{eqnarray}
where, $E$ is still central. The algebra (III.2) admit the following Hopf 
structure 
\begin{equation}
\begin{array}{lll}
\bigtriangleup_{h,w}(N)=N\otimes  e^{w A^{+}}+ 1\otimes N,&\qquad 
S_{h,w}(N)=-Ne^{-w A^{+}}, & \qquad \varepsilon_{h,w} (N)= 0, 
\nonumber  \\
\bigtriangleup_{h,w}(E)=E\otimes 1 + 1\otimes E, &\qquad S_{h,w}(E)=-E,
& \qquad \varepsilon_{h,w} (E)= 0, \nonumber \\
\bigtriangleup_{h,w} (A^{+})=A^{+} \otimes 1 + 1 \otimes A^{+}, &\qquad 
S_{h,w}(A^{+})=-A^{+}, & \qquad \varepsilon_{h,w} (A^{+})= 0,
\end{array}
\end{equation} 
and 
\begin{eqnarray}
&& \bigtriangleup_{h,w} (A)=A \otimes e^{hE}e^{w A^{+}}+e^{-hE}\otimes A
+we^{-hE}\;N\otimes {\sinh (hE) \over h}e^{w A^{+}},
\qquad\varepsilon_{h,w} (A)= 0,
\nonumber \\
&& \qquad\qquad \qquad \qquad S_{h,w}(A)=-Ae^{-w A^{+}}+wN\;{\sinh (hE) 
\over h} \;e^{-w A^{+}}.
\end{eqnarray}
\end{prop}
\noindent {\bf Proof:} All the Hopf algebra axioms can be verified by direct 
calculations. The elements 
$E$ and $A^{+}$ are primitives.\qed

\paragraph{}
The Heisenberg subalgebra generated by $E$, $A^{+}$ and $A$ is 
not a Hopf subalgebra. The Casimir element of the quantum algebra 
${\cal U}_{h,w}({\cal H}(4))$ is given by   
\begin{eqnarray}
C_{h,w}= N\; {\sinh{ h E} \over h} + {e^{-w A_{+}} -1 \over 2\;w} A+ 
A {e^{-w A_{+}} -1 \over 2\;w}.  
\end{eqnarray}
When $h$ is equal to zero, ${\cal C}_{h,w}$ correspond to the central element 
given in \cite{Her}. The explicit expression of $N$ in terms of $A$, $A^{+}$ 
and $E$ is given by the series 
\begin{eqnarray}
&&N\equiv\biggl({1- e^{-w A^{+}} \over w}\biggr)\;A\;\biggl[{\sinh (hE)\over h}
\biggr]^{-1} 
= 2h\biggl({1- e^{-w A^{+}} \over w}\biggr)
\;A\;\sum_{k=0}^{\infty}
e^{-(2k+1)h E}. 
\end{eqnarray}

The elements $E^{i}(A^{+})^j N^{k} A^{l}$ where $(i,j,k,l)\in\NN
\times \NN\times\NN\times \NN$ build a Poincar\'e-Birkhoff-Witt basis of 
${\cal U}_{h,w}({\cal H}(4))$:
  
\begin{prop}
Each infinite dimensional irreducible representation $\pi_{e,n}$ of 
${\cal U}_{h,w}({\cal H}(4))$ is labeled by two parameters $e$ and 
$n$. A generic representation $\pi_{e,n}$ is defined as follows:
\begin{eqnarray}
&&A|r\rangle=\biggl({\sinh(he)\over h}\biggr)^{1/2}\sqrt{r}|r-1\rangle, 
\nonumber\\
&& A^{+}|r\rangle=\sum_{k=0}^{\infty}{w^{k}\over (k+1)}
\biggl({\sinh(he)\over h}\biggr)^{(k+1)/2}
\sqrt{(r+k+1)!\over r!}|r+k+1\rangle, \nonumber\\
&& \biggl({1-e^{-wA^{+}}\over w}\biggr)|r\rangle =
\biggl({\sinh(he)\over h}\biggr)^{1/2}\sqrt{r+1}|r+1\rangle, \nonumber\\
&& E |r\rangle =e |r\rangle, \phantom{a \over b} \nonumber\\
&& N|r\rangle =(r+n) |r\rangle\phantom{a \over b}
\end{eqnarray}
where $\lbrace |r\rangle \rbrace_{r=0}^{\infty}$ is a orthonormal basis of the
module $V_{e,n}$.
\end{prop}

The invertible element ${\cal F}=e^{w N \otimes A^{+}}$ satisfies 
the relations (II.12)-(II.15) on the enveloping algebra 
${\cal U}_{h}({\cal H}(4))$. 
Now, Let us turn to the universal ${\cal R}$-matrix of ${\cal U}_{h,w}
({\cal H}(4))$:  
\begin{prop} 
The Hopf algebra ${\cal U}_{h,w}({\cal H}(4))$ is quasitriangular. The  
universal ${\cal R}$-matrix has the following form 
\begin{eqnarray} 
{\cal R}^{h,w}= (\sigma \circ {\cal F})^{-1}.{\cal R}^{h}. {\cal F}
\end{eqnarray} 
where 
\begin{eqnarray} 
&& {\cal F}= e^{w N \otimes A^{+}},\\ 
&& \displaystyle {\cal R}^{h}=
e^{-2h\;E\otimes N}\;\exp \biggl(2h\;e^{hE}A
\otimes \biggl({1-e^{-w A^{+}}\over w}\biggr)\biggr). 
\end{eqnarray}
The universal ${\cal R}$-matrix (III.8) solves the quantum Yang-Baxter 
equation (II.11).
\end{prop}
\noindent {\bf Proof:} The relations (II.9) are verified using the 
Campbell-Baker-Hausdorff formula 
\begin{eqnarray}
\displaystyle e^{\phi} \bigtriangleup(.)\;e^{-\phi}= \bigtriangleup(.) + 
\sum_{n=1}^{\infty} {1\over n!} \displaystyle \underbrace{[\phi,\cdots\;[\phi, 
\;\bigtriangleup(.)\;]\cdots ]}_{\hbox{n brackets}}.   
\end{eqnarray}
Recall that ${\cal F}$ satisfies the relations (II.12)-(II.15). 
Let us check that $(\bigtriangleup_{h,w}\otimes \id)({\cal R}^{h,w})=
{\cal R}^{h,w}_{13}{\cal R}^{h,w}_{23}$. This equation reduces to 
\begin{eqnarray}
&& {\cal F}_{31}^{-1}{\cal R}^{h}_{13}{\cal F}_{13}{\cal F}_{23}^{-1}
{\cal R}^{h}_{23}
{\cal F}_{23}={\cal F}_{12}^{-1}(\bigtriangleup_{h}\otimes \id)
({\cal F}_{21}^{-1})(\bigtriangleup_{h}\otimes \id)({\cal R}^{h})
(\bigtriangleup_{h}\otimes \id)({\cal F}){\cal F}_{12}, \nonumber\\ 
&& \phantom{{\cal F}_{31}^{-1}R^{h}_{13}{\cal F}_{13}{\cal F}_{23}^{-1}
{\cal R}^{h}_{23}{\cal F}_{23}}={\cal F}_{12}^{-1}(\bigtriangleup_{h}
\otimes \id)({\cal F}_{21}^{-1}){\cal R}^{h}_{13}{\cal R}^{h}_{23}
(\bigtriangleup_{h}\otimes \id)({\cal F}){\cal F}_{12}, \nonumber\\
&&\phantom{{\cal F}_{31}^{-1}R^{h}_{13}{\cal F}_{13}{\cal F}_{23}^{-1}
{\cal R}^{h}_{23}{\cal F}_{23}}={\cal F}_{12}^{-1}(\bigtriangleup_{h}
\otimes \id)({\cal F}_{21}^{-1}){\cal R}^{h}_{13}{\cal R}^{h}_{23}
(\id \otimes \bigtriangleup_{h})({\cal F}){\cal F}_{23}, \nonumber
\end{eqnarray}
namely,
\begin{eqnarray}
{\cal F}_{12}^{-1}(\bigtriangleup_{h}\otimes \id)({\cal F}_{21}^{-1})
{\cal R}^{h}_{13}= {\cal F}_{31}{\cal R}^{h}_{13}{\cal F}_{13}
{\cal F}_{23}^{-1}(\id \otimes \bigtriangleup_{h}')({\cal F}^{-1})  
\end{eqnarray}
where, $\bigtriangleup_{h}'=\sigma \circ\bigtriangleup_{h}$. Applying the 
flip $\sigma_{23}$ to both sides, the preceding equation is equivalent to
\begin{eqnarray}
&& {\cal F}_{13}^{-1}\sigma_{23}(\bigtriangleup_{h}\otimes\id)
({\cal F}_{21}^{-1}){\cal R}^{h}_{12}= {\cal F}_{21}{\cal R}^{h}_{12}
{\cal F}_{12}{\cal F}_{23}^{-1}(\id \otimes \bigtriangleup_{h})
({\cal F}^{-1}) \nonumber \\
&& \phantom{{\cal F}_{13}^{-1}\sigma_{23}(\bigtriangleup_{h}\otimes\id)
({\cal F}_{21}^{-1}){\cal R}^{h}_{12}}={\cal F}_{21}^{-1}{\cal R}_{12}^{h}
(\id \otimes \bigtriangleup_{h})({\cal F}^{-1}) \nonumber    
\end{eqnarray}  
namely,
\begin{eqnarray}
&& {\cal F}_{13}^{-1}\sigma_{23}(\bigtriangleup_{h}\otimes\id)
({\cal F}_{21}^{-1})= F_{21}(\bigtriangleup_{h}'\otimes \id)
({\cal F}^{-1}). 
\end{eqnarray}
It is easy to see that applying $\sigma_{12}$ to the left-hand side of the 
preceding equation gives ${\cal F}_{23}^{-1}(\id \otimes \bigtriangleup_{h})
({\cal F}^{-1})$, and applying it to the right-hand side gives 
${\cal F}_{12}^{-1}(\bigtriangleup_{h}\otimes \id)({\cal F}^{-1})$. So the 
result follows on using (II.13) once again. The relation  $(\id \otimes
\bigtriangleup_{h,w})({\cal R}^{h,w})={\cal R}^{h,w}_{13}{\cal R}^{h,w}_{23}$ 
is verified using the same method.   

\qed

The element (III.10) correspond to the ${\cal R}_{h}$-matrix (II.22) of the 
algebra (II.4) generated by $N$, $E$, $A$ and $(1-e^{-wA^{+}})/w$ (the 
generator $A^{+}$ is replaced by $(1-e^{-wA^{+}})/w$) , whereas 
${\cal F}^{-1}$ is the twist which deforms the coproducts (II.5) to the 
coproducts (III.3)-(III.4), namely 
\begin{eqnarray}
&& \bigtriangleup_{h,w}={\cal F}^{-1}\bigtriangleup_{h}{\cal F}.
\end{eqnarray}
Furthermore, 
\begin{eqnarray}
&& \varepsilon_{h,w}(\;.\;)=\varepsilon_{h}(\;.\;),\qquad\qquad 
 S_{h,w}(\;.\;)= v^{-1} S_{h}(\;.\;) v
\end{eqnarray}       
where, $v$ is the invertible element given by
\begin{eqnarray}
&&v=\sum_{k=0}^{\infty}{(-1)^{k}w^{k}\over k!}N^{k}(A^{+})^{k} \qquad\qquad
\nonumber \\ && v^{-1}=\sum_{k=0}^{\infty}{1\over k!}
N^{k}(\ln(2-e^{-wA^{+}}))^{k}  
\end{eqnarray}
which satisfies the following relation 
\begin{eqnarray}
v^{-1}S_{h}(v)=e^{wA^{+}}.
\end{eqnarray}

\paragraph{}
From the universal ${\cal R}$-matrix, by standard techniques, we can readily 
deduce a new ${\cal R}$-matrix depending on a continuous parameter $u$. In 
fact, defining the operator $T_{u}$ by its action: 
\begin{equation}
\begin{array}{lll}
\displaystyle T_{u}A^{+}=e^{-u}A^{+}, &\qquad  T_{u}w=e^{u}w, 
&\qquad T_{u}h=h, \phantom{a\over { c\over d}}\nonumber \\
T_{u}A=e^{u}A^{+}, &\qquad  T_{u}E=E, & \qquad T_{u}N=N,
\end{array}
\end{equation}
we can define
\begin{eqnarray}
{\cal R}^{h,w}(u)= (T_{u}\otimes 1){\cal R}^{h,w}= e^{-wA^{+}\otimes N}
e^{-2hE\otimes N}\exp \biggl(2he^{u}e^{hE}A
\otimes \biggl({1-e^{-w A^{+}}\over w}\Biggr)\biggr)
e^{wN\otimes A^{+}}\phantom{xxx}
\end{eqnarray}
and again by direct calculations we have: The matrix ${\cal R}^{h,w}(u)$ 
defined in (III.19) satisfies the Yang-Baxter equation 
\begin{eqnarray}
{\cal R}^{h,w}_{12}(u){\cal R}^{h,w}_{13}(u+v){\cal R}^{h,w}_{23}(v)=
{\cal R}^{h,w}_{23}(v){\cal R}^{h,w}_{13}(u+v){\cal R}^{h,w}_{12}(u). 
\end{eqnarray}
The classical $r(u)$-matrix corresponding to the universal 
${\cal R}^{h,w}(u)$-matrix
(III.19), depending on the parameter $u$, is the following 
\begin{eqnarray}
r(u)= w (n\otimes a^{+}-a^{+}\otimes n)+ 2\;h (e^{u}\; a\otimes a^{+}
-e\otimes n) 
\end{eqnarray}
and solves the parameter-dependent classical Yang-Baxter equation
 \begin{eqnarray}  
[r_{12}(u),r_{13}(u+v)]+[r_{12}(u),r_{23}(v)]+[r_{13}(u+v),r_{23}(v)]=0.  
\end{eqnarray}

\paragraph{}
Now, consider the infinite dimensional irreducible representations 
$\pi_{e,n}$ with the condition $e^{2he}\not=1$. Evaluates the interwiner 
matrices $R(e_1,e_2)\equiv e^{2he_1n_2}\sigma(\pi_{e_1,n_1}\otimes
\pi_{e_2,n_2}){\cal R}^{h,w}:\;V_{e_1,n_1}\otimes V_{e_2,n_2}\rightarrow 
V_{e_2,n_2}\otimes V_{e_1,n_1}$ and $R^{-1}(e_1,e_2)\equiv e^{-2he_2n_1}
(\pi_{e_1,n_1}\otimes\pi_{e_2,n_2})[{\cal R}^{h,w}]^{-1}\sigma:\;V_{e_1,n_1}
\otimes V_{e_2,n_2}\rightarrow V_{e_2,n_2}\otimes V_{e_1,n_1}$, we obtain  
\begin{eqnarray}
&&[R]_{r_1,r_2}^{r_1',r_2'}(e_1,e_2)=w^{r'_1+r'_2-r_1-r_2} 
\biggl({e^{he_1}- e^{-he_1}\over 2h}\biggr)^{(r'_2-r_1)/2}
\biggl({e^{he_2}- e^{-he_2}\over 2h}\biggr)^{(r'_1-r_2)/2}
\biggl({r'_1!r'_2!\over r_1!r_2!}\biggr)^{1/2} \nonumber \\
&&\phantom{(R^{h,w})_{r_1,r_2}^{r_1',r_2'}=} 
\sum_{s=\hbox{sup}(r_1-r'_2,0)}^{\hbox{inf}(r_1,r'_1-r_2)}
\biggl({r_1 \atop s}\biggr)\phi_{r'_1-r_2-s}{\bar\phi}_{r'_2-r_1+s}
(e^{he_1}- e^{-he_1})^s e^{h(s-2r_1')e_1} , \\
&&[R^{-1}]_{r_1,r_2}^{r_1',r_2'}(e_1,e_2)=w^{r'_1+r'_2-r_1-r_2}
\biggl({e^{he_1}-e^{-he_1}\over 2h}\biggr)^{(r'_2-r_1)/2}
\biggl({e^{he_2}-e^{-he_2}\over 2h}\biggr)^{(r'_1-r_2)/2}
\biggl({r'_1!r'_2!\over r_1!r_2!}\biggr)^{1/2} \nonumber \\
&&\phantom{(R^{h,w})_{r_1,r_2}^{r_1',r_2'}=}
\sum_{s=\hbox{sup}(r_2-r'_1,0)}^{\hbox{inf}(r_2,r'_2-r_1)}
(-1)^s\biggl({ r'_1+s \atop s}\biggr)\phi_{r'_1-r_2+s}
{\bar\phi}_{r'_2-r_1-s}(e^{he_2}- e^{-he_2})^{s}e^{h(s+2 r_1)e_2},
\end{eqnarray}
with
\begin{eqnarray}
&&\phi_{r'_1-r_2\pm s}=\sum_{k=0}^{r'_1-r_2\pm s}{f_{k}^{r'_1-r_2\pm s-k} 
\over k!}(r_1+n_1)^{k},\qquad 
{\bar\phi}_{r'_2-r_1\pm s}=\sum_{k=0}^{r'_2-r_1\pm s}{(-1)^{k}
f_{k}^{r'_2-r_1\pm s-k}\over k!}(r'_1+n_2)^{k},\qquad \nonumber 
\end{eqnarray}
where 
\begin{eqnarray}
f_{k}^{s}=\displaystyle\left\{\matrix{ 1 && \hbox {if}\qquad k=0\cr 
&& \cr 
\displaystyle \sum_{i_1+i_2+\cdots + i_k=s} 
{1\over (i_1+1)(i_2+1)\cdots (i_k+1)} &&\hbox {if}\qquad
k \geq 1 \cr}\right. 
\end{eqnarray}
Taking $w=0$, the relations (III.23) and (III.24) reduces to the results 
obtained by Gomez and Sierra in \cite{Gom}, namely
\begin{eqnarray}
&& [R_g]_{r_1,r_2}^{r_1',r_2'}(e_1,e_2)=\delta_{r'_1+r'_2,r_1+r_2} 
\biggl({r_1 \atop r_2'}\biggr)^{1/2}
\biggl({r'_1 \atop r_2}\biggr)^{1/2}e^{-h(r_1'+r_2)e_1} \nonumber \\
&& \phantom{[R_g]_{r_1,r_2}^{r_1',r_2'}(e_1,e_2)=} 
\biggl[(e^{he_1}- e^{-he_1})(e^{he_2}-e^{-he_2})
\biggr ]^{(r'_1-r_2)/2}  \nonumber \\
&& [R_g^{-1}]_{r_1,r_2}^{r_1',r_2'}(e_1,e_2)=\delta_{r'_1+r'_2,r_1+r_2} 
\biggl({r_2 \atop r_1'}\biggr)^{1/2}\biggl({r'_2 \atop r_1}\biggr)^{1/2}
(-1)^{r_2-r'_1}e^{h(r_2'+r_1)e_2} \nonumber \\
&& \phantom{ [R_g]_{r_1,r_2}^{r_1',r_2'}(e_1,e_2)=} 
\biggl[(e^{he_1}- e^{-he_1})(e^{he_2}- e^{-he_2})\biggr ]^{(r_2-r'_1)/2}.
\end{eqnarray}

\paragraph{}
Recall that, a braid group ${\cal B}_{m}$ is an abstract group generated by 
the elements $\sigma_i,\;1\leq i\leq m-1$ satisfying the relations $\sigma_i
\sigma_j=\sigma_j\sigma_i$ if $|i-j |\geq 2$, $\sigma_{i}\sigma_{i}^{-1}=
\sigma_{i}^{-1}\sigma_{i}=1$ and $\sigma_i\sigma_{i+1}\sigma_i=\sigma_{i+1}
\sigma_i \sigma_{i+1}$ if $1\leq i\leq m-2$. In colored braid groups cases, 
$\sigma_i$ contain different nontrivial parameters called the string variables.
The $R$ matrix (III.23) build a new infinite dimensional representation 
$\rho_{m}$ of the colored braided group ${\cal B}_{m}$ via
\begin{eqnarray}
&&\rho_{m}:{\cal B}_{m}\longrightarrow\hbox{End}(\otimes_{i=1}^{m}
V_{e_i,n_i}) \nonumber \\
&&\phantom{\rho_{m}:}\;\;\sigma_{i} \;\longmapsto  \id^{\otimes(i-1)}
\otimes R^{e_i,e_{i+1}}\otimes\id^{\otimes(m-i-1)} 
\end{eqnarray} 
which satisfies
\begin{eqnarray}
R_{i}^{e_2,e_3}R_{i+1}^{e_1,e_3}R_{i}^{e_1,e_2}=
R_{i+1}^{e_1,e_2}R_{i}^{e_1,e_3}R_{i+1}^{e_2,e_3}. 
\end{eqnarray} 
The noncolored version correspond to take $e_1=e_2=\cdots = e_m=e$.  

Now, let us concentrate in noncolored case. The noncolored braid group 
representation (NCBG) admits an extension \`a la Turaev \cite{Tur}, if there 
exists an isomorphism $\mu: V_{e,n}\rightarrow V_{e,n}$ which transforms 
the basis $\{|r\rangle\}_{r=0}^{\infty}$ into $\{\mu_{r}|r\rangle\}_{
r=0}^{\infty}$ satisfying the following three conditions:   
\begin{eqnarray}
&& (\mu_i \mu_j-\mu_k \mu_l) R_{i,j}^{k,l}(e,e)=0, \phantom{\sum_{j}}
\qquad\qquad\qquad \phantom{a}\\
&& \sum_{j} R_{i,j}^{k,j}(e,e)\mu_j = \delta_{i}^{k} ab  ,  \\
&&  \sum_{j} [R^{-1}]_{i,j}^{k,j}(e,e)\mu_j = \delta_{i}^{k} a^{-1}b 
\end{eqnarray} 
where $a,b$ are constants. The Turaev conditions (III.29)-(III.31) hold if 
\begin{eqnarray}
\mu=\id,\qquad\qquad a=e^{he}\qquad\qquad b= e^{-he}.
\end{eqnarray}
Let $\zeta:{\cal B}_{m}\rightarrow\ZZ$ such that $\zeta(\sigma_i^{\pm 1})=
\pm 1$, $\zeta(xy)=\zeta(x)+\zeta(y)$. Then, the link invariant
$P:\prod_{m\geq 2}{\cal B}_{m}\rightarrow \CC$ associated to the enhanced 
Yang-Baxter operator (EYB-operator) $(R, \id, e^{he}, e^{-he})$ is 
\begin{eqnarray}
P(x)=q^{e(-\zeta(x)+m)}Tr[\rho_{m}(x)],\qquad\qquad\forall x\in{\cal B}_{m}
\end{eqnarray}
which is invariant under the two Markov moves $P(xy)=P(yx)$, $x,y\in 
{\cal B}_{m}$ (type I) and $P(x\sigma_{m}^{\pm 1})=P(x),\;x\in{\cal B}_{m},
\;\sigma_{m}\in {\cal B}_{m+1}$ (type II). Let us remark that: 

\smallskip
\smallskip

{\bf (i)} The combination of the standard and nonstandard quantizations
permit to build a new infinite dimensional representation 
of the colored braided group ${\cal B}_{m}$.

\smallskip
\smallskip

{\bf (ii)} The link invariant (III.33) is equal to the one obtained in 
reference \cite{Gom}. The link invariants obtained for a structure 
defind as combination of the standard and nonstandard quantizations 
are exactly equal to those calculated for the standard quantization 
separetely. The invariants (III.33) are polynomials only in the variable 
$q^{he}$. Its easy to see that the nonstandard parameter $w$ have no 
contribution in (III.33) (see the results (III.23) and (III.24)).  

\smallskip
\smallskip

{\bf (iii)} These comments are true for the colored braiding version. 
(See reference \cite{Gom} for more details).

\sect{The Two Parametric Heisenberg Group $H_{h,w}(4)$}  

A $3\times 3$ dimensional representation $\pi_{3}$ of the two parametric 
deformed Heisenberg algebra ${\cal U}_{h,w}({\cal H}(4))$ is given by
\begin{eqnarray}
&& \pi_{3}(A)=\pmatrix{0 & 1& 0\cr
 0 & 0& 0\cr 0 & 0& 0\cr}, \qquad\qquad \pi_{3}(A^{+})=\pmatrix{0 & 0& 0\cr
 0 & 0& 1\cr 0 & 0& 0\cr},   \nonumber \\
&& \pi_{3}(E)=\pmatrix{0 & 0& 1\cr
 0 & 0& 0\cr 0 & 0& 0\cr}, \qquad\qquad \pi_{3}(N)=\pmatrix{0 & 0& 0\cr
 0 & 1& 0\cr 0 & 0& 0\cr}.
\end{eqnarray}
This represenation remains undeformed. Correspondingly, the $R$-matrix 
(III.8) is represented by the $9\times 9$ matrix
\begin{eqnarray}
  && R^{h,w}_{\pi_{3}}=
(\pi_{3}\otimes \pi_{3})({\cal R}^{h,w})=\pmatrix{I_{3}&&2h\pi_{3}(A^{+})& 
&-2h\pi_{3}(N)\cr
 && && \cr
 0 &&  I_{3}+\omega\pi_{3}(A^{+}) &&  -w\pi_{3}(N) \cr 
 && && \cr
0 && 0&& I_{3}\cr},
\end{eqnarray}
$I_{3}$ being the $3\times 3$ identity matrix. We shall now present the 
quantum group $H(4)$ as a matrix quantum group \`a la Woronowicz 
\cite{Wor1,Wor2}. Let us consider a $3\times 3$ matrix of the following form 
\begin{eqnarray}
&&T=e^{\beta\otimes \pi_{3}(E)}e^{\delta\otimes\pi_{3}(A^{+})}
e^{\gamma\otimes\pi_{3}(N)} e^{\alpha\otimes\pi_{3}(A)}=
\pmatrix{1 & \alpha & \beta\cr
0 & e^{\gamma} & \delta \cr 
0 & 0 & 1 \cr},
\end{eqnarray}
where the matrix elements $\alpha$, $\beta$, $\gamma$ and $\delta$ generate 
the algebra of functions on the quantum group ${\cal F}_{h,w}(H(4))$. The 
matrix elements $\alpha$, $\beta$, $\gamma$ and $\delta$ of $T$ satisfies 
the relation: 
\begin{eqnarray} 
R^{h,w}_{\pi_{3}}\;T_{1}\;T_{2}=T_{2} \;T_{1}\;R ^{h,w}_{\pi_{3}}
\end{eqnarray}
where, $T_{1}=T\otimes 1$ and $T_{2}=1\otimes T$. The algebra ${\cal F}_{h,w}
(H(4))$ can be endowed with a Hopf algebra by defining a comultiplication 
$\bigtriangleup$ and counit $\varepsilon$ as 
\begin{eqnarray}
\bigtriangleup(T)=T{\dot\otimes}T,\qquad\qquad\qquad\varepsilon(T)=1,   
\end{eqnarray}
where ${\dot\otimes}$ denotes matrix multiplication and tensor product 
of the $\CC^{\star}$-algebras of noncommutative representative functions. 
The inverse matrix then defines the antipode, namely
\begin{eqnarray} 
S(T)=T^{-1}.  
\end{eqnarray}

\begin{prop} \hfill 

(i) The relations between the generators $\alpha$, $\beta$, $\gamma$ and 
$\delta$ are the following:
\begin{eqnarray}
&&[\alpha,\beta]=2h\alpha+w\alpha^2,\qquad\qquad\qquad [\gamma,\delta]=
w(e^{\gamma}-1),\nonumber\\
&& [\alpha,\delta]= w \alpha e^{\gamma},\qquad\;\;\qquad\qquad\qquad 
[\beta, \gamma]= - w \alpha , \nonumber \\
&& [\beta,\delta]=[\alpha,\gamma ]=0. 
\end{eqnarray}

(ii) The generators $\alpha$, $\beta$, $\gamma$ and $\delta$ with the 
relations specified in (IV.7) have coproducts, counits and antipodes given by
\begin{equation}
\begin{array}{lll}   
\bigtriangleup (\alpha) = \alpha \otimes e^{\gamma} +1\otimes \alpha, &\qquad
S(\alpha)= - e^{-\gamma} \;\alpha,&\qquad \varepsilon (\alpha)=0, \nonumber \\
\bigtriangleup(\beta)=\beta\otimes 1+1\otimes\beta+\alpha\otimes\delta,&\qquad
S(\beta)=-\beta+e^{-\gamma}\alpha\delta,&\qquad\varepsilon(\beta)=0,\nonumber\\
\bigtriangleup (e^{\gamma}) = e^{\gamma} \otimes e^{\gamma}, &\qquad
S(e^{\gamma})=e^{-\gamma} , & \qquad \varepsilon (\gamma)=0,\nonumber \\
\bigtriangleup(\delta)=\delta\otimes 1+e^{\gamma}\otimes\delta, &\qquad
S(\delta)= - e^{-\gamma} \delta, &\qquad \varepsilon (\delta)=0.
\end{array}
\end{equation}
\end{prop}

\noindent {\bf Proof:} The relations (IV.7) arises from the relation (IV.4). 
The Hopf structure of ${\cal F}_{h,w}(H(4))$ reads from the equations (IV.5) 
and (IV.6), namely 
\begin{eqnarray}
&& T^{-1}=\pmatrix{1 & -e^{-\gamma}\alpha & -\beta+e^{-\gamma}
\alpha \delta \cr
0 &  e^{-\gamma} &  -e^{-\gamma}\delta \cr 
0 & 0& 1\cr}
\end{eqnarray}
and
\begin{eqnarray}
&& T{\dot \otimes} T = 
\pmatrix{1\otimes 1  & \alpha \otimes e^{\gamma} +1\otimes \alpha & 
 \beta\otimes 1+1\otimes\beta+\alpha\otimes\delta  \cr
0&e^{\gamma}\otimes e^{\gamma}&\delta\otimes 1+e^{\gamma}\otimes\delta  \cr 
0 & 0& 1\otimes 1\cr}.
\end{eqnarray}
\qed

When the parameter $h$ is equal to zero, ${\cal F}_{h,w}(H(4))$ is reduced 
to the algebra ${\cal F}_{w}(H(4))$ obtained by Ballesteros and Herranz 
\cite{Her}:
\begin{eqnarray}
&&[\alpha,\beta]=w\alpha^2,\qquad\;\qquad\qquad [\gamma, \delta ]=
w(e^{\gamma}-1),\nonumber\\
&& [\alpha,\delta]=w\alpha e^{\gamma},\qquad\qquad \qquad
[\beta ,\gamma ]=-w\alpha, \nonumber \\
&& [\alpha ,\gamma ]=[\beta,\delta ]=0. 
\end{eqnarray}  
Whereas, if $w=0$, the algebra ${\cal F}_{h}(H(4))$ being the standard ones 
(according to the algebra obtained in \cite{Cel}, we take $\gamma=0$), i.e.   
 \begin{eqnarray}
&&[\alpha,\beta]=2h\alpha,\qquad\qquad\qquad 
[\alpha ,\delta ]=[\beta ,\delta ]=0
\end{eqnarray}
and 
\begin{equation}
\begin{array}{lll}   
\bigtriangleup (\alpha) = \alpha \otimes 1 +1\otimes \alpha, &\qquad
S(\alpha)= -  \alpha,&\qquad \varepsilon (\alpha)=0, \nonumber \\
\bigtriangleup(\beta)=\beta\otimes 1+1\otimes\beta+\alpha\otimes\delta,&\qquad
S(\beta)=-\beta+\alpha \delta,&\qquad\varepsilon(\beta)=0,\nonumber \\
\bigtriangleup(\delta)=\delta\otimes 1+1\otimes\delta, &\qquad
S(\delta)= -  \delta, &\qquad \varepsilon (\delta)=0.
\end{array}
\end{equation}
The relations (IV.13) are obtained in the reference \cite{Cel}. The algebra 
(IV.12) is lightly different from the algebra ${\cal F}_{h}(H(4))$ obtained 
in \cite{Cel}.

\sect{The Ribbon Hopf Algebra ${\cal U}_{h,w}({\cal H}(4))$}

Any quasitriangular Hopf algebra (${\cal U}$, ${\cal R}$) has an invertible 
element, usually called $u$, with the property that
\begin{eqnarray}
S^{2}(a)= uau^{-1}, \qquad\qquad \forall a \in {\cal U}. 
\end{eqnarray}   
The element $u$ and its inverse $u^{-1}$ can be obtained from the universal 
${\cal R}$-matrix as follows:
\begin{eqnarray} 
&& u = m (S\otimes \id)({\cal R}_{21}), \nonumber\\
&& u^{-1}= m(S^{-1} \otimes \id)({\cal R}_{21}^{-1}). 
\end{eqnarray}
The elements $u$ and $uS(u)$ satisfies the relations \cite{Turaev}
\begin{eqnarray}
&&\varepsilon(u)=1,\qquad\qquad \bigtriangleup(u)=({\cal R}_{21}{\cal R})^{-1}
(u\otimes u)=(u\otimes u) ({\cal R}_{21}{\cal R})^{-1}, \nonumber \\
&& \bigtriangleup(uS(u))=({\cal R}_{21}{\cal R})^{-2}
(u\otimes u). 
\end{eqnarray}
The element $uS(u)=S(u)u$ is central in ${\cal U}$. The elements $u$ and
$S(u)$ of ${\cal U}$ does not in general commute with the others elements 
of ${\cal U}$, however $u$ and $S(u)$ are central if $S^{2}=\id$.

Recall that, a Ribbon Hopf algebra $({\cal U}, {\cal R}, \theta)$ is a 
quasitriangular Hopf algebra with a central element $\theta$ satisfying
\begin{eqnarray}
&&\theta^{2}=u S(u),\qquad\qquad\qquad\qquad\qquad
\qquad S(\theta)=\theta,\nonumber\\
&&\bigtriangleup(\theta)=({\cal R}_{21}{\cal R})^{-1}(\theta\otimes\theta), 
\qquad\qquad\qquad  \varepsilon(\theta)=1. 
\end{eqnarray}

\begin{prop}
Let (${\cal U}$, ${\cal R}$) be a quasitriangular Hopf algebra and ${\cal F}$ 
an invertible element of ${\cal U}\otimes {\cal U}$ satisfying (II.12)-(II.15). 
The $u_{\cal F}$ 
element of the twisted quasitriangular Hopf algebra ${\cal U}_{\cal F}$ is 
\begin{eqnarray}
&& u_{\cal F}= v^{-1}\;S(v)\;u=v^{-1}uS^{-1}(v), 
\end{eqnarray}    
where $v$ and its inverse are given by (II.17). The 
element $u_{\cal F}$ has the property that
\begin{eqnarray}
&&S^{2}_{\cal F}(a)=u_{\cal F}\;a\;u_{\cal F}^{-1},\qquad\qquad a\in 
{\cal U}_{\cal F},\\
&& \epsilon_{\cal F}(u_{\cal F})=1, \\
&& \bigtriangleup_{\cal F}(u_{\cal F})= ({\cal R}^{\cal F}_{21}
{\cal R}^{\cal F})^{-1}(u_{\cal F}\otimes u_{\cal F})=(u_{\cal F}\otimes 
u_{\cal F})({\cal R}^{\cal F}_{21}{\cal R}^{\cal F})^{-1}, 
\end{eqnarray}   
\end{prop}  

\noindent {\bf Proof:} For the proof, we choose ${\cal F}=\sum_{(f)}
f_1\otimes f_2$, ${\cal F}^{-1}=\sum_{(g)}g_1\otimes g_2$, ${\cal R}=
\sum_{(r)}r_1\otimes r_2$, $u=\sum_{(r)}S(r_2)r_1$, $v=\sum_{(f)}S(f_1)f_2$
and $v^{-1}=\sum_{(g)}g_1S(g_2)$. The element $u_{\cal F}$ is obtained from 
the universal ${\cal R}_{\cal F}$-matrix as follows:
\begin{eqnarray}
&& u_{\cal F}= m(S_{\cal F}\otimes \id)({\cal R}^{\cal F}_{21}), \nonumber \\
&&\phantom{u_{\cal F}}= m(v^{-1}\otimes\id)(S\otimes \id)({\cal F}^{-1}
{\cal R}_{21}{\cal F}_{21})(1\otimes v), \nonumber \\
&&\phantom{u_{\cal F}}= \sum_{(f)}\sum_{(g)}\sum_{(r)} 
v^{-1}S(f_2)S(r_2)S(g_1) v g_2 r_1 f_1. \nonumber
\end{eqnarray}
From the relation $m(S\otimes \id)({\cal F}{\cal F}^{-1})= 1$, we obtain 
that $\sum_{(g)}S(g_1) v g_2=1$. Now, the preceding relations reads     
\begin{eqnarray}
&& u_{\cal F}= \sum_{(f)}\sum_{(r)} v^{-1}S(f_2)u f_1\nonumber \\
&& \phantom{u_{\cal F}}=\sum_{(f)}\sum_{(r)} 
v^{-1}S(f_2) S^{2}(f_1)u. 
\end{eqnarray}
So the result follows on using the element $v$ (II.17). The equations (V.6) 
and (V.7) are verified by direct calculations. The relations (V.8) arise 
as follows: 
\begin{eqnarray}
&& \bigtriangleup_{\cal F}(u_{\cal F})={\cal F}^{-1}
\bigtriangleup(v^{-1}S(v)u)\;{\cal F}, \nonumber \\
&& \phantom{\bigtriangleup_{\cal F}(u_{\cal F})}= {\cal F}^{-1}
\bigtriangleup(v^{-1}S(v))\bigtriangleup(u)\;{\cal F}, \nonumber \\
&& \phantom{\bigtriangleup_{\cal F}(u_{\cal F})}=
(v^{-1}S(v)\otimes v^{-1}S(v))((S^{2}\otimes S^{2}){\cal F}^{-1})
\bigtriangleup(u)\;{\cal F}, \nonumber \\
&& \phantom{\bigtriangleup_{\cal F}(u_{\cal F})}=
(v^{-1}S(v)u\otimes v^{-1}S(v)u){\cal F}^{-1}
({\cal R}_{21}{\cal R})^{-1}{\cal F},\nonumber \\
&& \bigtriangleup_{\cal F}(u_{\cal F})=(u_{\cal F}\otimes 
u_{\cal F})({\cal R}^{\cal F}_{21}{\cal R}^{\cal F})^{-1}. 
\end{eqnarray}
\qed

\begin{prop} 
Let $({\cal U}, {\cal R}, \theta)$ be a ribbon Hopf algebra and ${\cal F}$ an
invertible element of ${\cal U} \otimes {\cal U}$ satisfying (II.12)-(II.15). 
The twisted algebra ${\cal U}_{\cal F}$ is a ribbon Hopf algebra with 
$\theta_{\cal F}=\theta$ and the relations
\begin{eqnarray}
&& {\cal S}_{\cal F}(\theta_{\cal F})=1,\qquad\qquad\qquad\qquad\qquad\qquad
\varepsilon_{\cal F}(\theta_{\cal F})=1, \nonumber \\
&& \bigtriangleup_{\cal F}(\theta_{\cal F})= ({\cal R}^{\cal F}_{21}
{\cal R}^{\cal F})^{-1}(\theta_{\cal F}\otimes \theta_{\cal F}). 
\end{eqnarray}
\end{prop}

\noindent {\bf Proof:} The results (V.11) arises from $\theta^{2}_{\cal F}=
u_{\cal F}S_{\cal F}(u_{\cal F})= uS(u)$, namely, $\theta_{\cal F}=\theta$.    

\qed

In the particular case of the universal ${\cal R}_{h}$-matrix (II.22) of the 
Heisenberg algebra ${\cal U}_{h}({\cal H}(4))$ we obtain from proposition 
9 and \cite{Gom}  
\begin{eqnarray}
&& u_{h}=\sum_{l\geq 0}{(-1)^{l} (2h)^{l} \over l!} 
e^{-hlE} (A^{+})^{l} A^{l} e^{2hEN}, \nonumber\\
&& u_{h}^{-1}= \sum_{l\geq 0}{(2h)^{l} \over l!} 
e^{hlE} (A^{+})^{l} A^{l} e^{-2hEN}, 
\end{eqnarray}
For the standard Heisenberg algebra ${\cal U}_{h}({\cal H}(4))$, we have from 
(II.5) that $S^{2}=\id$ so that $u_h$ is central. Similarly, $S_{h}(u_h)$ is 
also central and we have 
\begin{eqnarray} 
S_{h}(u_h)=e^{-2hE}u_h. 
\end{eqnarray} 
Gomez and Sierra have proved that ${\cal U}_{h}({\cal H}(4))$ is a Ribbon 
Hopf algebra with $\theta_{h} = e^{-hE}u_h$ \cite{Gom}. For the two parametric 
Heisenberg algebra ${\cal U}_{h,w}({\cal H}(4))$, the $u_{h,w}$ element and
its inverse reads 
\begin{eqnarray}
&& u_{h,w}=e^{wA^{+}} \sum_{l\geq 0}{(-1)^{l} (2h)^{l} \over l!} 
e^{-hlE}\biggl({1-e^{-wA^{+}}\over w}\biggr)^{l} A^{l} 
e^{2hEN}, \nonumber\\
&& u_{h,w}^{-1}=\sum_{l\geq 0}{(2h)^{l} \over l!} 
e^{hlE} \biggl({1-e^{-wA^{+}}\over w}\biggr)^{l} A^{l} e^{-2hEN}
e^{-wA^{+}},
\end{eqnarray}
and 
\begin{eqnarray}
S_{h,w}(u_{h,w})= e^{-2hE}e^{-wA^{+}} u_{h,w}.
\end{eqnarray}
For the representation $\pi_{e,n}$, we obtain 
\begin{eqnarray}
&& u_{h,w}|r\rangle =e^{2nhe}\sum_{l\geq 0}w^{l}
\biggl({t-t^{-1}\over 2h}\biggr)^{l/2}
\biggl({(r+l)! \over r!}\biggr)^{1/2} |r+l\rangle , \nonumber \\
&& u_{h,w}^{-1}|r\rangle =e^{-2nhe}|r\rangle -wt^{-2n} 
\biggl({\sinh(he)\over h}\biggr)^{1/2} \sqrt{r+1}|r+1\rangle.
\end{eqnarray}

\begin{prop}
The quasitriangular Hopf algebra $({\cal U}_{h,w}({\cal H}(4)),
{\cal R}^{h,w})$ is a Ribbon Hopf algebra with
\begin{eqnarray}
&& \theta_{h,w}= e^{-hE} \sum_{l\geq 0}{(-1)^{l} (2h)^{l} \over l!} 
e^{-hlE}\biggl({1-e^{-wA^{+}}\over w}\biggr)^{l}A^{l}e^{2hEN}\nonumber\\
&& \theta_{h,w}^{-1}= e^{hE} \sum_{l\geq 0}{(2h)^{l} \over l!} 
e^{hlE}\biggl({1-e^{-wA^{+}}\over w}\biggr)^{l}A^{l}e^{-2hEN}
\end{eqnarray}
which in the irreducible representation $\pi_{e,n}$ takes the value 
\begin{eqnarray}
\theta_{h,w}|r\rangle =e^{(2n-1)he}|r\rangle, \qquad\qquad 
\theta_{h,w}^{-1}|r\rangle =e^{-(2n-1)he}|r\rangle. 
\end{eqnarray}
\end{prop}

Finally, let us remark that: 

\smallskip
\smallskip

{\bf (i)} The eigenvalue of the central element $\theta_{h,w}$ and its inverse 
depend only on the standard parameter $h$.

\smallskip
\smallskip

{\bf (ii)} The value of $\theta_{h,w}$ in a given irreducible representation 
$\pi_{e,n}$ contains some interesting information of the corresponding 
conformal field theory (CFT) associated to the quantum algebra ${\cal U}_{h,w}
({\cal H}(4))$. In the cases of a semisimple algebra ${\cal G}$ it was shown 
that the conformal weight $\Delta_{\alpha}$ of a primary field of the WZW 
model ${\hat {\cal G}}_{k}$ is related to the value $\theta_{\alpha}$ by 
\cite{Gom1}
\begin{eqnarray}
\theta_{\alpha}= e^{2\pi i \Delta_{\alpha}} 
\end{eqnarray}      
where $\theta_{\alpha}$ is the value of $\theta$ on the irreducible 
representation $\alpha$ of ${\cal U}_{q_1,\cdots, q_{\xi}}({\cal G})$ 
associated to the primary field $\alpha$. If the formula (V.19) holds 
true for the quantum Heisenberg algebra ${\cal U}_{h,w}({\cal H}(4))$ 
it would imply that 
\begin{eqnarray}
e^{(2n-1)he}= e^{2\pi i \Delta_{e,n}}. 
\end{eqnarray}    
The conformal weight $\Delta_{e,n}$ depend {\em only} on the standard 
parameter $h$. The nonstandard parameter $w$ has no contribution 
in WZW model.

\sect{Conclusions and Perspectives}

Let us remark that the combination of the standard and nonstandard 
deformations in the case of the enveloping Heisenberg algebra  is 
possible because the invertible element ${\cal F}$ satisfies the 
axioms (II.12)-(II.15) on the standard algebra ${\cal U}_{h}({\cal H}(4))$. 
The parameter arising from the nonstandard quantization not appear 
in the conformal weight $\Delta_{e,n}$ of the primary field of the WZW 
model ${\hat {\cal G}}_{k}$ and in the link invariants. The parameters
which plays a relevant roles arises from the standard quantization.          
In the case of the quantum algebra ${\cal U}_{q}(sl(2))$, the analogs
two parametric deformation cannot be defined because the twist ${\cal D}$
do not obeys to the axioms (II.12)-(II.15).   
        
Using the annihilator and creator operators of the infinite Heisenberg 
algebra ${\cal U}_{h,w}({\cal H}(\infty))$ and the Sugawara construction,
the two parametric deformation can be extended to Virasoro algebra. The two 
parametric deformation of the Galilei group obtained by contraction can be 
also used to study the magnetic chain following the approach developed 
in \cite{Bon}. These problems will be studied else where.

\vskip 2cm 
\noindent {\bf Acknowledgments:} 
 
The author wants to thank Daniel Arnaudon, Amitabha Chakrabarti and 
Chryssomalis Chryssomalakos for interesting discussions. 
He is also grateful to Ranabir Chakrabarti for a precious help.

\vskip 1cm

\newpage

\bibliographystyle{amsplain}

\end{document}